\documentclass{article}
\usepackage{amssymb}
\usepackage{graphics}
\input xy
\xyoption{all}

\begin{document}
\title{Classification of $1^{st}$ order symplectic spinor operators over   contact projective geometries}
\author{Svatopluk Kr\'ysl \footnote{{\it E-mail address}: krysl@karlin.mff.cuni.cz}\\ {\it \small  Charles University, Sokolovsk\'a 83, Praha, Czech Republic}
\\{\it \small and
Humboldt-Universit\"{a}t zu Berlin, Unter den Linden 6, Berlin, Germany.}
\thanks{I am very grateful to Vladim\' \i r Sou\v{c}ek for useful comments, Andreas \v{C}ap for introducing me into
contact projective geometries and to David Vogan who recommended me
some of his texts on globalization of Harish-Chandra modules for
reading. The author of this article was supported by the grant
GA\v{C}R 201/06/P223 for young researchers of The Grant Agency of
Czech Republic and by the grant GA UK 447/2004. Supported also by the SPP 1096 of the DFG.}}

\maketitle

\noindent\centerline{\large\bf Abstract}

We give a classification of $1^{st}$ order invariant differential operators
acting between sections of certain bundles associated to Cartan geometries of
the so called metaplectic contact projective type.
 These bundles are associated via representations, which
 are derived from the so called  higher symplectic, harmonic or
 generalized Kostant spinor modules.  Higher symplectic spinor
 modules
are arising from the Segal-Shale-Weil representation of the
metaplectic group by tensoring it by finite dimensional modules. We
show that for all pairs of the considered bundles, there is at most
one $1^{st}$ order invariant differential operator up to a complex
multiple and give an equivalence condition for the existence of such
an operator. Contact projective analogues of the well known Dirac,
twistor and Rarita-Schwinger operators appearing in Riemannian
geometry are special examples of these operators.

{\it Math. Subj. Class.:} 22E46, 58J60, 58J70.


\section{Introduction}

   The operators we would like to classify are   $1^{st}$ order invariant differential operators  acting between sections of vector bundles
associated to metaplectic contact projective geometries via certain minimal globalizations.

    Metaplectic contact projective geometry on an odd dimensional manifold is   first   a contact geometry,
    i.e., it is given by a corank one subbundle of the tangent bundle of the manifold which is
     nonintegrable in the Frobenius sense in each point of the manifold. Second part of the
     metaplectic contact projective structure on a manifold   is given by a   class  of projectively
      equivalent contact  partial affine connections. Here, partial contact means that the connections
       are compatible with the contact structure and that they are acting only on the sections of the
       contact subbundle. These connections are called projectively equivalent because they have the same
       class of unparameterized geodesics going in the contact subbundle direction,  see, e.g.,
       D. Fox \cite{Fox}, where you can find a relationship between the contact projective geometries and
        classical path geometries. The adjective "metaplectic" suggests that in addition
    to contact projective geometries, the metaplectic contact
        projective structures    include some spin phenomena like the spin structures over Riemannian manifolds.
Metaplectic contact projective and contact projective geometries
have their description also via Cartan geometries. Contact
projective geometries could be modeled on a $(2l+1)$-dimensional
projective space $\mathbb{P}\mathbb{V}$ of a $(2l+2)$-dimensional
real symplectic vector space $\mathbb{V},$ which we suppose to be
equipped with a symplectic form $\omega.$ Here, the projective space
is considered as a homogeneous space $G/P,$ where $G$ is the
symplectic Lie group $Sp(\mathbb{V},\omega)$  acting transitively on
$\mathbb{P}\mathbb{V}$ by the factorization of its defining
representation (on $\mathbb{V}),$ and $P$ is an isotropy subgroup of
this action.  In this case, it is easy to see that $P$ is a
parabolic subgroup, which turns out to be crucial for our
classification. Contact projective geometry, in the sense of \'{E}.
Cartan, are curved versions $(p: \mathcal{G}\to M,\omega)$ of this
homogeneous (also called Klein) model $G/P$.
 There exist certain conditions (known as normalization conditions) under which the Cartan's
 principal bundle approach and  the classical one (via the class of connections and the contact subbundle)
 are equivalent, see, e.g., \v{C}ap, Schichl \cite{Cap} for details. We also remind that contact geometries
  are an arena for time-dependent Hamiltonian mechanics.  Klein model of the metaplectic contact
  projective geometry consists of two groups $\tilde{G}$ and $\tilde{P},$ where $\tilde{G}$ is the
   metaplectic group $Mp(\mathbb{V},\omega),$ i.e., a non-trivial double covering of the symplectic group $G,$ and $\tilde{P}$ is the preimage of $P$ by this covering.

Symplectic spinor operators over projective contact geometries are
acting between sections of the so called higher symplectic spinor
bundles.
     These bundles are associated via certain infinite dimensional irreducible admissible representations of the parabolic principal group $P.$
    The   parabolic group $P$ acts then non trivially only by  its Levi
factor $G_0$, while the action of the unipotent part is trivial. The
semisimple part $\mathfrak{g}_0^{ss}$ of the Lie algebra  of the
Levi part of the parabolic group $P$ is isomorphic to the symplectic
Lie algebra $\mathfrak{sp}(2l,\mathbb{R}).$ Thus to give an
admissible representation of $P,$ we have to specify a
representation of $\mathfrak{g}_0^{ss}.$ Let us recall that the
classification of first order invariant operator was done by
Slov\'{a}k, Sou\v{c}ek in \cite{SS} (generalizing an approach of
Fegan in \cite{Fegan}) for all finite dimensional irreducible
representations  and general parabolic subgroup $P$ of a semisimple $G$ (almost Hermitian structures are studied in
detail).
  Nevertheless, there are some interesting infinite dimensional representations of the complex symplectic Lie algebra, to which we shall focus our attention. These
representations form a class consisting of infinite dimensional modules with bounded
multiplicities. Modules with bounded
multiplicities are representations, for which
there is a nonnegative integer,  such that the dimension of each
weight space of this module is bounded by it from above. Britten,
Hooper and Lemire in \cite{BHL} and Britten, Hooper in \cite{BH}
showed that each of these modules appear as direct summands in a
tensor product of a finite dimensional $\mathfrak{sp}(2l,\mathbb{C})$-module and the so called Kostant (or basic) symplectic spinor
module $\mathbb{S}_+$ and vice versa.  Irreducible representations
in this completely reducible tensor product are called   higher
symplectic,  harmonic   or  generalized Kostant spinors. It is well known, that all
finite dimensional modules over complex symplectic Lie algebra
appear as  irreducible submodules of a tensor power of the defining
representation.
 Thus the infinite dimensional modules with bounded multiplicities are analogous to
  the spinor-vector representations of complex orthogonal Lie algebras.
 Namely, each finite dimensional module over orthogonal Lie algebra is an
irreducible summand in the tensor product of a basic spinor
representation and some power of the defining module (spinor-vector representations),
or in the power of the defining
representation itself (vector representations). In order to
have a complete picture, it remains to show that the
   basic (or Kostant) spinors  are analogous to the
orthogonal ones, even though infinite dimensional. The basic symplectic
spinor module $\mathbb{S}_+$ was discovered by Bertram Kostant, when
he was introducing half-forms  for metaplectic structures over
symplectic manifolds in the context of geometric quantization. While
in the orthogonal case spinor representations can be realized using the
exterior algebra of a maximal isotropic vector space, the symplectic spinor
representations are realized  using the  symmetric algebra of certain maximal isotropic vector space (called Lagrangian in the symplectic setting). This procedure goes
roughly as follows: one takes the Chevalley realization of the
symplectic Lie algebra $C_l$ by polynomial coefficients linear
differential operators acting on polynomials $\mathbb{C}[z^1,\ldots,
z^l]$ in $l$ complex variables. The space of polynomials splits into
two irreducible summands over the symplectic Lie algebra, namely
into the two basic symplectic spinor modules $\mathbb{S}_+$ and
$\mathbb{S}_-$.  There is a relationship between the modules
$\mathbb{S}_+$ and $\mathbb{S}_-$ and the Segal-Shale-Weil or
oscillator representation. Namely, the underlying $C_l$-structure of
the Segal-Shale-Weil representation is isomorphic to
$\mathbb{S}_+\oplus \mathbb{S}_-.$

  In order to classify  $1^{st}$ order invariant differential operators, one needs to understand the structure of the space of $P$-homomorphisms between the so called $1^{st}$ jets prolongation  $P$-module of the domain module
and the target representation of $P$, see chapter 4. Thus the
classification problem translates into an algebraic one. In our
case, representation theory teaches us, that it is sometimes
sufficient to understand our representation at its infinitesimal
level. The only thing one needs in this case, is to understand the
infinitesimal version of the $1^{st}$-jets prolongation module. For
our aims, the most important part    of the $1^{st}$ jets
prolongation  module  consists of a tensor product of the defining
representation of $C_l$ and a higher symplectic spinor module. In
order to describe the space of $P$-homomorphisms, one needs to
decompose the mentioned tensor product into irreducible summands.
This was done by Kr\'{y}sl in \cite{KryslJOLT}, where results of
Humphreys in \cite{Hum2} and Kac and Wakimoto in \cite{KW} were
used.

   Let us mention that some of these operators   are contact analogues of the well known
 symplectic Dirac operator, symplectic Rarita-Schwinger and symplectic twistor operator. Analytical
   properties of these operators were studied by many authors, see, e.g., K. Habermann \cite{KH} and A.
   Klein \cite{Klein}. These symplectic versions were mentioned also by M. B. Green and C. M. Hull,
   see \cite{GH}, in the context of covariant quantization of 10 dimensional super-strings and also in
    the theory of Dirac-K\"{a}hler fields, see Reuter \cite{Reuter}, where we
    found a motivation for our studies of this topic.

In the second section, metaplectic contact projective geometries are
defined using the Cartan's approach. Basic properties of higher
symplectic spinor modules (Theorem 1) together with a theorem on a
decomposition of the tensor product of the defining representation
of $\mathfrak{sp}(2l,\mathbb{C})$ and an arbitrary higher symplectic
spinor module (Theorem 2) are summarized in section 3 . Section 4 is
devoted to the classification result. Theorem 3 and Lemmas 1 and 2
in this section are straightforward generalizations of similar
results obtained by Slov\'{a}k and Sou\v{c}ek in \cite{SS}. Theorem
4 (in section 4) is a well known theorem on the action of a Casimir
element on highest weight modules. While in the subsection 4.1. we
are interested only in the classification at the infinitesimal level
(Theorem 5), we present our classification theorem at the globalized
level in subsection 4.2 (Theorem 6). In the fifth section, three
main examples of the $1^{st}$ order symplectic spinor operators over
contact projective structures are introduced.

\section{Metaplectic contact projective geometry}

The aim of this section is neither to serve as a comprehensive introduction into metaplectic  contact projective geometries, nor to
list all references related to this subject. We shall only present  a definition of metaplectic  contact projective geometry
by introducing its Klein model, and give only a few references, where one can find links to a broader literature on this topic (contact projective geometries, path geometries e.t.c.).

For a fixed positive integer $l \geq 3,$ let us consider a real
symplectic vector space $(\mathbb{V}, \omega)$ of real dimension
$2l+2$  together with the defining action of the symplectic Lie
group $G:=Sp(\mathbb{V},\omega).$  The defining action is transitive
on $\mathbb{V}-\{0\},$ and thus it defines a transitive action $G
\times \mathbb{PV} \to \mathbb{PV}$ on the projective space
$\mathbb{PV}$ of $\mathbb{V}$  by the prescription $(g,[v]) \mapsto
[gv]$ for $g\in G$ and $v \in \mathbb{V}-\{0\}.$ (Here, $[v]$
denotes the one dimensional vector subspace  spanned by $v.$) Let us
denote the stabilizer of a point in $\mathbb{PV}$ by $P.$ It is well
known that this group is a parabolic subgroup of $G,$ see, e.g., D.
Fox \cite{Fox}. The pair $(G,P)$ is often called Klein pair of
contact projective geometry. Let us denote the Lie algebra of $P$ by
$\mathfrak{p}.$

{\bf Definition 1:} Cartan geometry $(p: \mathcal{G} \to M^{2l+1}, \omega)$ is called a contact projective geometry of rank $l$, if
it is a Cartan geometry modeled on the Klein geometry of type $(G,P)$ for $G$ and $P$ introduced above.

It is possible to show that each contact projective geometry defines
 a contact structure on the tangent bundle $TM$ of the base manifold $M$ and a class $[\nabla]$ of
contact projectively   equivalent partial affine
connections $\nabla$ acting on the sections of the contact subbundle
(see the Introduction for some remarks).
 For more details on this topic, see Fox \cite{Fox}. In \v{C}ap, Schichl
 \cite{Cap}, one can find a treatment on the equivalence problem for contact projective structures.
   Roughly speaking, the reader can find a proof there, that under certain conditions, there is an  isomorphism
 between the Cartan approach and the classical one (via contact subbundle and a class of connections).
Because we would like to include some spin phenomena, let us consider a slightly modified situation.
Fix a non-trivial two-fold covering $q: \tilde{G} \to G$ of the symplectic group $G=Sp(\mathbb{V},\omega)$ by the metaplectic group  $\tilde{G}=Mp(\mathbb{V},\omega),$ see Kashiwara, Vergne \cite{KV}. Let us denote the $q$-preimage of $P$ by $\tilde{P}.$

{\bf Definition 2:} Cartan geometry $(p:\tilde{\mathcal{G}}\to M^{2l+1},\omega)$ is called metaplectic contact projective geometry of rank $l$, if it is a Cartan geometry modeled
on the Klein geometry of type $(\tilde{G},\tilde{P})$ with $\tilde{G}$ and $\tilde{P}$ introduced above.

Let us remark, that in definition 2,
we do not demand     metaplectic contact projective structure to be connected  to a contact projective structure as one demands in the case of spin structures over Rimannian manifolds or in the case of metaplectic structures over manifolds with a symplectic structure.

\section{Higher symplectic spinor modules}
  Let $C_l \simeq \mathfrak{sp}(2l,\mathbb{C}),$ $l \geq 3,$
  be the complex   symplectic Lie algebra. Consider a Cartan subalgebra $\mathfrak{h}$ of $C_l$  together with a choice of positive roots
$\Phi^+.$  The set of fundamental
weights $\{\varpi_i\}_{i=1}^l$ is then uniquely determined. For
later use, we shall need an orthogonal basis (with respect to the
  form dual to the Killing form of $C_l$), $\{\epsilon_i\}_{i=1}^l,$
for which $\varpi_i=\sum_{j=1}^i\epsilon_j$ for $i=1,\ldots, l.$

For $\lambda \in \mathfrak{h}^*,$ let $L(\lambda)$ be the  irreducible $C_l$-module
with the highest  weight  $\lambda.$ This module is defined uniquely up to a $C_l$-isomorphism.
If $\lambda$ happens to be integral
and dominant (with respect to the choice of $(\mathfrak{h},\Phi^+)$),
i.e., if  $L(\lambda)$ is finite dimensional, we shall write
$F(\lambda)$ instead of $L(\lambda).$ Let $L$ be an arbitrary
(finite or infinite dimensional) weight module over a complex simple
Lie algebra.
 We call $L$ a module with bounded multiplicities, if there is a $k\in \mathbb{N}_0,$ such that for each $\mu \in \mathfrak{h}^*,$ $\mbox{dim}\,L_{\mu}\leq k,$ where $L_{\mu}$ is the weight space of weight $\mu.$

Let us introduce the following set of weights
$$A:=\{\lambda=\sum_{i=1}^{l}\lambda_i\varpi_i|\lambda_i \in \mathbb{N}_0, i=1,\ldots, l-1, \lambda_{l-1}+2\lambda_l + 3 >0,
 \lambda_l \in \mathbb{Z}+\frac{1}{2}\}.$$

{\bf Definition 3:} For a weight $\lambda \in
A,$  we call the  module $L(\lambda)$ higher symplectic spinor module. We shall denote the module
$L(-\frac{1}{2}\varpi_{l})$ by $\mathbb{S}_+$   and the module $L(\varpi_{l-1}-\frac{3}{2}\varpi_l)$ by
$\mathbb{S}_-.$ We shall call these two representations basic
symplectic spinor modules.

The next theorem says that the class of higher symplectic spinor
modules is quite natural and   in a sense broad.

{\bf Theorem 1:} Let   $\lambda \in \mathfrak{h}^*$. Then the following are equivalent
\begin{itemize}
\item[1.)] $L(\lambda)$ is an infinite dimensional $C_l$-module with bounded multiplicities;
\item[2.)] $L(\lambda)$ is a direct summand  in $\mathbb{S}_+\otimes F(\nu)$ for some integral dominant $\nu \in \mathfrak{h}^*;$
\item[3.)] $\lambda \in A.$
\end{itemize}
{\it Proof.} See Britten, Hooper, Lemire \cite{BHL} and Britten, Lemire \cite{BH}. $\Box$

In the next theorem, the tensor product of a higher symplectic
spinor module and the defining representation $\mathbb{C}^{2l}
\simeq F(\varpi_1)$ of the complex symplectic Lie algebra $C_l$ is
decomposed into irreducible summands. We shall need this statement
in the classification procedure. It gives us an important
information on the structure of the $1^{st}$ jets prolongation
module for metaplectic contact projective structures.

{\bf Theorem 2:} Let $\lambda \in A.$ Then
$$L(\lambda) \otimes F(\varpi_1)=\bigoplus_{\mu \in A_{\lambda}}L(\mu),$$
where $A_{\lambda}:=A \cap \{\lambda + \nu| \nu \in \Pi(\varpi_1)
\}$ and $\Pi(\varpi_1)=\{\pm\epsilon_i| i=1,\ldots,l\}$ is the   set of weights of the
defining representation.

{\it Proof.} See Kr\'{y}sl,  \cite{KryslJOLT}. $\Box$

Let us remark, that the proof of this theorem is based on the so called Kac-Wakimoto formal character formula   published in \cite{KW} (generalizing
a statement of Jantzen in \cite{Jan}) and   some results of Humphreys, see \cite{Hum2}, in which results of Kostant (from \cite{K}) on tensor products of finite and infinite dimensional modules admitting a central character are specified.

\section{Classification of first order invariant operators}

In this section, we will be investigating first order invariant differential
operators acting between  sections of certain vector bundles associated to parabolic
geometries $(p:\mathcal{G} \to M,\omega),$ i.e., to Cartan geometries modeled on Klein
 pairs $(G,P),$ where $P$ is an arbitrary parabolic subgroup of an arbitrary semisimple
 Lie group $G.$

We first consider a general real semisimple Lie group $G$ together with its parabolic subgroup $P$
 and then we restrict our attention to
  the metaplectic contact projective case.
Let us suppose that    the Lie algebra $\mathfrak{g}$   of the group
$G$ is equipped with a $|k|$-grading
$\mathfrak{g}=\bigoplus_{i=-k}^k\mathfrak{g}_i,$ i.e.,
$\mathfrak{g}_1$ generates $\bigoplus_{i=1}^k\mathfrak{g}_i$ as a
Lie algebra and $[\mathfrak{g}_i,\mathfrak{g}_j] \subseteq
\mathfrak{g}_{i+j}$ for $i,j \in \{-k,\ldots, k\}.$ \footnote{By
definition, $\mathfrak{g}_i=0$ for $|i|>k$ is to be understood.}
Denote the semisimple part and the center of the reductive Lie
algebra $\mathfrak{g}_0\subset \mathfrak{g}$ by
$\mathfrak{g}_0^{ss}$ and  $\mathfrak{z}(\mathfrak{g}_0),$
respectively. The subalgebra $\bigoplus_{i=0}^k\mathfrak{g}_i$ forms
a parabolic subalgebra of $\mathfrak{g}$ and will be denoted by
$\mathfrak{p}.$ Let us suppose that $\mathfrak{p}$ is isomorphic to
the      Lie algebra of the fixed parabolic subgroup $P$ of $G.$ The
nilpotent part $\bigoplus_{i=1}^k\mathfrak{g}_i$ of $\mathfrak{p}$
is usually denoted by $\mathfrak{g}_+$ and the negative $\bigoplus_{i=-k}^{-1}\mathfrak{g}_i$
part of $\mathfrak{g}$  by $\mathfrak{g}_-$. Let us
consider   Killing forms $(,)_{\mathfrak{g}}$ and
$(,)_{\mathfrak{g}_0^{ss}}$ of $\mathfrak{g}$ and
$\mathfrak{g}_0^{ss},$ respectively. Further, fix a basis
$\{\xi^i\}_{i=1}^{r}$ of $\mathfrak{g}_{+},$ such that
$\{\xi^{i}\}_{i=1}^s$ is a basis of $\mathfrak{g}_1$ and
$\{\xi^{i}\}_{i=s+1}^r$ is a basis of
$\bigoplus_{i=2}^k\mathfrak{g}_i.$ The second basis, we will use, is
a basis of $\mathfrak{g}_0^{ss},$ which will be denoted by
$\{\eta^i\}_{i=1}^t.$ The $|k|$-grading of $\mathfrak{g}$ uniquely
determines the so called grading element $Gr \in
\mathfrak{z}(\mathfrak{g}_0).$ The defining equation for this
element is $[Gr,X] = jX$ for $X \in \mathfrak{g}_j$ and each $j \in
\{-k,\ldots,k\}.$ It is known that for each $|k|$-grading of a real (or complex) semisimple Lie algebra the grading element exists, see, e.g., Yamaguchi \cite{Yama}.
Sometimes, we will denote the grading element $Gr$
by $\eta^{t+1}.$ The set $\{\eta^i\}_{i=1}^{t+1}$ is then a basis of
$\mathfrak{g}_0.$ Let us denote the basis of $\mathfrak{g}_-$ dual
to $\{\xi^i\}_{i=1}^{r}$ with respect to the Killing form
$(,)_{\mathfrak{g}}$ by $\{\xi_i\}_{i=1}^{r}$ and the basis of
$\mathfrak{g}_0$ dual to the basis $\{\eta^i\}_{i=1}^{t+1}$ with
respect to the Killing form $(,)_{\mathfrak{g}}$ by
$\{\eta_{i}\}_{i=1}^{t+1}.$

At the beginning, let us consider two complex irreducible
representations $(\sigma, \bf E)$ and $(\tau, \bf F)$ of $P$
  in the category $\mathcal{R}(P),$ the objects of which are locally convex, Hausdorff vector spaces with a continuous
linear action of $P,$ which is admissible, of finite length.
Here, admissible action means that the restriction of this action to the Levi subgroup $G_0$ of $P$ is admissible, see Vogan \cite{Vogan}.
The
morphisms in the category $\mathcal{R}(P)$ are linear continuous
$P$-equivariant maps between the objects. It is well known that the
unipotent part of the parabolic group acts trivially on both $\bf E$
and $\bf F.$ We shall call ${\bf E}$ and ${\bf F}$ the domain and
the target module, respectively and we shall specify further
conditions on these representations later.
 Generally, for a Lie group $G$ and its admissible representation
${\bf E},$ we shall denote the corresponding Harish-Chandra $(\mathfrak{g},K)$-module ($K$ is maximal compact in $G$) by $E$  and when we will only  be considering the $\mathfrak{g}$-module structure, we shall use the symbol $\mathbb{E}$ for it.
Further, we will denote the corresponding actions of an element $X$ from the Lie algebra of $G$ on a vector $v$ simply by $X.v,$ and the action of $g \in G$ on a vector $v$ by $g.v$ - the considered representation will be clear from a context.

Let us stress that most our proofs are formally
almost identical to that ones written by Slov\'{a}k, Sou\v{c}ek in \cite{SS}, but we formulate them also for infinite dimensional admissible irreducible $\bf{E}$ and $\bf{F}$, and use the decomposition result in Kr\'{y}sl \cite{KryslJOLT} when we will be treating the metaplectic contact projective case.

Let $(p:\mathcal{G} \to M, \omega)$ be a Cartan geometry modeled on
the Klein pair $(G, P).$ Because $\omega_u: T_{u}\mathcal{G} \to
\mathfrak{g}$ is an isomorphism for each $u \in \mathcal{G}$   by
definition, we can define a vector field $\omega^{-1}(X)$ for each
$X \in \mathfrak{g}$ by the equation
$\omega_{u}(\omega^{-1}(X)_u)=X,$ the so called constant vector
field. For later use, consider two associated vector bundles ${\bf
E}M:=\mathcal{G}\times_{\sigma}{\bf E}$ and ${\bf
F}M:=\mathcal{G}\times_{\tau}{\bf F}$ - the so called domain and
target bundle, respectively. To each Cartan geometry, there is an
associated  derivative $\nabla^{\omega}$ defined as follows. For any section
$s \in \Gamma(M, {\bf E}M)$ considered as $s \in
\mathcal{C}^{\infty}(\mathcal{G},{\bf E})^P$ under the obvious
isomorphism, we obtain a mapping $\nabla^{\omega}s: \mathcal{G} \to
\mathfrak{g}_-^*\otimes {\bf E},$ defined by the formula
$$(\nabla^{\omega}s(u))X:=\mathcal{L}_{\omega^{-1}(X)}s(u),$$ where
$X \in \mathfrak{g}_-,$  $u \in \mathcal{G}$ and $\mathcal{L}$ is
the Lie derivative. The associated derivative $\nabla^{\omega}$ is
usually called absolute invariant derivative.
The $1^{st}$ jets prolongation module $J^1{\bf E}$ of ${\bf E}$ is
defined as follows. As a  vector space, it is simply the space ${\bf
E}\oplus (\mathfrak{g_+}\otimes {\bf E}).$ To be specific, let us
fix the Grothendieck's projective tensor product topology on
$1^{st}$ jets prolongation module, see Treves \cite{Tr} or/and D.
Vogan \cite{Vogan}. The vector space $J^1{\bf E}$ comes up with an
inherited natural action of the group $P,$  forming the $1^{st}$
jets prolongation $P$-module,  see \v{C}ap, Slov\'{a}k, Sou\v{c}ek
\cite{SlovakLectures}. Let us remark that the function
$u\mapsto(s(u),\nabla^{\omega}s(u))$ defines a $P$-equivariant
function on $\mathcal{G}$ with values in $J^1{\bf E}$ and thus a
section of the first jet prolongation bundle $J^1({\bf E}M)$ of the
associated bundle ${\bf E}M.$ For details, see \v{C}ap, Slov\'{a}k
Sou\v{c}ek \cite{SlovakLectures}.

 By differentiation of the $P$-action on $J^1{\bf E}$, we can obtain a  $\mathfrak{p}$-module structure, the so called infinitesimal $1^{st}$ jets
  prolongation $\mathfrak{p}$-module $J^1\mathbb{E},$ which is as a vector space isomorphic
to $\mathbb{E}\oplus (\mathfrak{g}_+\otimes \mathbb{E}).$  The $\mathfrak{p}$-representation is then given by the formula
\begin{equation}
R.(v',S\otimes v''):=(R.v',S\otimes R.v''+[R,S]\otimes v''+\sum_{i=1}^r\xi^i\otimes [R,\xi_i]_{\mathfrak{p}}.v') \label{akce}
\end{equation}
where $R \in \mathfrak{p},$ $S\in \mathfrak{g}_+,$  $v', v'' \in
\mathbb{E}$ and $[R,\xi_i]_{\mathfrak{p}}$ denotes the projection of
$[R,\xi_i]$   to $\mathfrak{p}.$ For a derivation of the above
formula, see \v{C}ap, Slov\'{a}k, Sou\v{c}ek \cite{SlovakLectures}
for more details.
 Obviously, this action does not depend on a choice of the vector space basis $\{\xi^i\}_{i=1}^{r}.$ We will call this action
  the induced action of  $\mathfrak{p}.$

  {\bf Definition 4:} We call a vector space homomorphism
  $\mathfrak{D}: \Gamma(M,{\bf E} M) \to \Gamma(M, {\bf F} M)$
   first order invariant differential operator,
 if there is a $P$-module homomorphism\footnote{By a $P$-module homomorphism,
 we mean a morphism in $\mathcal{R}(P).$} ${\bf D}: J^1{\bf E} \to {\bf F}$,
 such that $\mathfrak{D}s(u)={\bf D}(s(u), \nabla^{\omega}s(u))$ for each
 $u \in \mathcal{G}$ and each section $s \in \Gamma(M, {\bf E}M)$
  (considered as a $P$-equivariant ${\bf E}$-valued smooth function on $\mathcal{G}).$

Let us remark, that this definition could be generalized for an arbitrary order. The  corresponding operators are called strongly invariant.
There exist also operators which are invariant in a broader sense (see \v{C}ap, Slov\'{a}k, Sou\v{c}ek \cite{CSS}) and not strongly invariant.

 We shall denote the vector space of first order invariant differential operators by $\hbox{Diff}({\bf E}M, {\bf F}M)^{1}_{(p:\mathcal{G}\to M)}.$ It is clear that $\hbox{Diff}({\bf E}M, {\bf F}M)^{1}_{(p:\mathcal{G}\to M)} \simeq \mbox{Hom}_{P}(J^1{\bf E},{\bf F})$ as complex vector spaces.
Let us denote the restricted $1^{st}$ jets prolongation $P$-module, i.e., the quotient $P$-module
$$[{\bf E}\oplus (\mathfrak{g}_+\otimes {\bf E})]/[\{0\}\oplus (\bigoplus_{i=2}^k\mathfrak{g}_i \otimes {\bf E})],$$
by $J^1_R {\bf E}.$
According to our notation, the meanings of $J_{R}^1 E$ and $J_{R}^1 \mathbb{E}$ are also fixed.
Now, let us introduce a linear mapping $\Psi: \mathfrak{g}_1\otimes \mathbb{E}\to \mathfrak{g}_1\otimes \mathbb{E}$
given by the following formula
$$\Psi(X\otimes v):=\sum_{i=1}^s\xi^i\otimes[X,\xi_i].v.$$
Obviously,  mapping $\Psi$ does not depend on a choice of the basis
$\{\xi^i\}_{i=1}^s$.

First, let us derive the following

{\bf Theorem 3:} Let $\mathbb{E}$ and $\mathbb{F}$ be two $\mathfrak{p}$-modules such that the nilpotent part $\mathfrak{g}_+$ acts trivially on them. If $D\in \mbox{Hom}_{\mathfrak{p}}(J^1\mathbb{E},\mathbb{F})$ is a $\mathfrak{p}$-homomorphism, then $D$ vanishes on the image of $\Psi$ and
$D$ factors through the restricted jets, i.e., $D(0 , Z \otimes v'')=0$ for each $v''  \in \mathbb{E}$ and $Z\in \bigoplus_{i=2}^k\mathfrak{g}_{i}.$
   Conversely, suppose $D\in \mbox{Hom}_{\mathfrak{g}_0}(J^1\mathbb{E},\mathbb{F})$ is a $\mathfrak{g}_0$-homomorphism,
$D$ factors through the restricted jets,
and $D$ vanishes on the image of $\Psi,$ then $D$ is a $\mathfrak{p}$-module homomorphism.

{\it Proof.}   Let $D\in
\mbox{Hom}_{\mathfrak{p}}(J^1\mathbb{E},\mathbb{F})$ be a
$\mathfrak{p}$-homomorphism. Take an element $\tilde{v}\in
\mathfrak{g}_+ . J^1\mathbb{E}.$ Then $D(\tilde{v})=D(X.v)$ for some
$X\in \mathfrak{g}_+$ and $v \in \mathbb{E}.$  Using the fact, that
$D$ is a $\mathfrak{p}$-homomorphism, we can write
$D(\tilde{v})=X.D(v)=0,$ because  the nilpotent algebra
$\mathfrak{g}_+$ acts trivially on the module $\mathbb{F}.$ Thus $D$
vanishes on the image of $\mathfrak{g}_+$ on $J^1\mathbb{E}.$

Now, we would like to prove, that $D$ factors through
$J^1_R\mathbb{E}$. Take an arbitrary element $Z\in
\bigoplus_{i=2}^k\mathfrak{g}_i$ and $v''\in \mathbb{E}.$ Because
$\mathfrak{g}$ is a $|k|$-graded algebra, there are $n \in \mathbb{N}$ and $X_i,Y_i \in
\mathfrak{g}_+$ for $i=1,\ldots, n,$ such that
$Z=\sum_{i=1}^{n}[X_i, Y_i].$ It is easy to compute that
$\sum_{i=1}^n X_i.(0,Y_i\otimes v'')=(0,\sum_{i=1}^n Y_i \otimes
X_i.v'' +[X_i,Y_i].v''+ 0)=(0,\sum_{i=1}^n[X_i,Y_i]\otimes
v'')=(0,Z\otimes v'').$ Thus we may write $D(0,Z\otimes
v'')=D(\sum_{i=1}^nX_i.(0,Y_i \otimes v''))=0,$ because $D$ acts
trivially on $\mathfrak{g}_+.J^1\mathbb{E},$ as we have already
proved.

Second, we shall prove that $D$ vanishes on the image of $\Psi.$
Substituting $v''=0$ into formula (\ref{akce}) for the induced action,
we get  that $X.(v',
0) = (X.v',\sum_{i=1}^r\xi^i\otimes [X,\xi_i]_{\mathfrak{p}}.v')$ for $v'\in \mathbb{E}$ and $X\in \mathfrak{g}_1.$
Assuming that the nilpotent subalgebra $\mathfrak{g}_+$ acts
trivially on $\mathbb{E}$, one obtains $X.(v',0) = (0,
\sum_{i=1}^r\xi^i\otimes [X,\xi_i]_{\mathfrak{p}}.v') = (0,
\sum_{i=1}^s\xi^i\otimes [X,
\xi_i]_{\mathfrak{p}}.v'+\sum_{i=s+1}^{r}\xi^{i}\otimes [X,
\xi_{i}]_{\mathfrak{p}}.v').$ The last summand is zero, because $[X,
\xi_i]_{\mathfrak{p}}=0$ for $i>s.$ Thus we have $X.(v',0) = (0,
\sum_{i=1}^s\xi^i\otimes [X,\xi_i]_{\mathfrak{p}}.v').$ Because $D$
vanishes on the image of the action of $\mathfrak{g}_+$ on
$J^1\mathbb{E},$ we know that $0 =D(X.(v',
0))=D(0,\sum_{i=1}^s\xi^i\otimes [X,\xi_i]_{\mathfrak{p}}.v').$
Since one can omit the restriction of the Lie bracket in the last
term to the subalgebra $\mathfrak{p}$ (we are considering
$\xi_{i}$ only for $i= 1,\ldots, s$), $D$ vanishes on the image of
$\Psi.$

Now, we would like to prove  the opposite direction. Hence suppose,
a $\mathfrak{g}_0$-homomorphism $D$ is given. Let us take an element
$S\in \mathfrak{g}_+$ (for $S\in \mathfrak{g}_0$ it is clear) and an
arbitrary element $\tilde{v}=(v',Y\otimes v'') \in J^1\mathbb{E}.$
Thus $D(S.\tilde{v})=D(S.(v',Y\otimes v''))=D(S.v',Y\otimes
S.v''+[S,Y]\otimes v''+ \sum_{i=1}^r\xi^i\otimes
[S,\xi_i]_{\mathfrak{p}}. v')= D(0, \sum_{i=1}^r\xi^i\otimes
[S,\xi_i]_{\mathfrak{p}}. v')=0=S.D(\tilde{v}),$ where we have used
that the action of $\mathfrak{g}_+$ is trivial on $\mathbb{E},$ $D$
factors through the restricted jets, vanishes on the image of
$\Psi,$  and the fact that $\mathfrak{g}_+$ acts trivially on
$\mathbb{F}$ . $\Box$

  Now, we derive the following

{\bf Lemma 1:} For the  mapping $\Psi,$ we have
 $$\Psi(X\otimes v)=\sum_{j=1}^{t+1} [\eta_j,X]\otimes
\eta^j.v$$ for each $X \in \mathfrak{g}_1$ and $v\in \mathbb{E}.$

{\it Proof.} Take an element $X \in \mathfrak{g}_1.$ Using the invariance of the Killing form $(,)_{\mathfrak{g}},$ expressed by
$[X,\xi_i]=\sum_{i=1}^{t+1}(\eta_i,[X,\xi_i])_{\mathfrak{g}}\eta^i=\sum_{i=1}^{t+1}([\eta_i,X],\xi_i)_{\mathfrak{g}}\eta^i,$ we compute the value $\Psi(X\otimes v)$ as
\begin{eqnarray*}
\Psi(X\otimes v) &=&\sum_{i=1}^s\xi^i\otimes
[X,\xi_i].v\\
&=&\sum_{i=1}^s\xi^i\otimes
\sum_j^{t+1}(\eta_j,[X,\xi_i])_{\mathfrak{g}}\eta^j.v\\
&=&\sum_{i=1}^s\xi^i\otimes\sum_{j=1}^{t+1}([\eta_j,X],\xi_i)_{\mathfrak{g}}\eta^j.v\\
&=&\sum_{i=1}^s\sum_{j=1}^{t+1}([\eta_j,X],\xi_i)_{\mathfrak{g}}\xi^i\otimes
\eta^j.v\\
&=& \sum_{j=1}^{t+1} [\eta_j, X] \otimes \eta^j.v.
\end{eqnarray*}
$\Box$

For any real Lie algebra $\mathfrak{g},$ let us denote its complexification over reals by
$\mathfrak{g}^{\mathbb{C}},$ i.e., $\mathfrak{g}^{\mathbb{C}}=\mathfrak{g}\otimes_{\mathbb{R}}\mathbb{C}.$
Let $\mathfrak{h}$ be a (complex) Cartan subalgebra of $(\mathfrak{g}_0^{ss})^{\mathbb{C}}.$
 For each $\lambda, \mu, \alpha \in \mathfrak{h}^{*},$ we  define a complex number
 $$c_{\lambda\alpha}^{\mu}=\frac{1}{2}[(\lambda,\lambda+2\delta)_{\mathfrak{g}_0^{ss}}+(\alpha,\alpha+2\delta)_{\mathfrak{g}_0^{ss}}
-(\mu,\mu+2\delta)_{\mathfrak{g}_0^{ss}}],$$ where $\delta$ denotes
the   sum of fundamental weights with respect to a choice of
positive roots. \footnote{We are denoting the Killing form on
$\mathfrak{g}_0^{ss}$ as well as the dual form on
$(\mathfrak{g}_0^{ss})^*$ by the same symbol
$(,)_{\mathfrak{g}_0^{ss}}.$ We shall also not distinguish between the Killing form of a real algebra and that one of the complexification of this algebra. We hope that this will cause no confusion.}

From now on, we shall suppose that the semisimple part of the Levi
factor of $P$ is actually
 simple and the center of the   Levi factor is one dimensional. These assumptions
  are rather technical and introduced only in order to simplify   formulations of our
  statements.
   Until yet, we have demanded the considered modules to be admissible irreducible
$P$-modules. In particular, we have used the fact that     the
unipotent part of $P$ acts trivially on them. From now on, we will suppose
in addition that the modules $\mathbb{E}$ and $\mathbb{F}$ are irreducible highest weight modules over     the complexification
  $(\mathfrak{g}_0^{ss})^{\mathbb{C}}$ of the Lie algebra $\mathfrak{g}_0^{ss}$ of the semisimple part of the Levi factor
  $G_0$ of $P.$   Further we shall suppose, that
  the grading element acts by a complex multiple on each of the modules $\mathbb{E}$ and $\mathbb{F}$.   We call a pair $(\lambda,
c) \in \mathfrak{h}^*\times \mathbb{C}$ a highest weight of a
representation ${\mathbb E}$ over the reductive Lie algebra
$(\mathfrak{g}_0)^{\mathbb{C}},$ if the restriction of the
representation of $\mathfrak{g}_0$ on $\mathbb{E}$ to the simple
part $(\mathfrak{g}_0^{ss})^{\mathbb{C}}$ has highest weight
$\lambda$ and the grading element $Gr$ acts by a complex number $c.$ The complex number $c$ is often called generalized
conformal weight of the $\mathfrak{p}$-module $\mathbb{E}.$

Recall a well known theorem on the action of the universal Casimir
element on highest weight modules.

{\bf Theorem 4:} Let
$\mathbb{E}$ be a highest weight module over the simple complex Lie
algebra $(\mathfrak{g}_0^{ss})^{\mathbb{C}}$ with a highest weight $\lambda \in
\mathfrak{h}^*$ and $C \in \mathfrak{U}((\mathfrak{g}_0^{ss})^{\mathbb{C}})$ be the
universal Casimir element of $(\mathfrak{g}_0^{ss})^{\mathbb{C}}.$ Then
$$C.v=(\lambda,\lambda+2\delta)_{\mathfrak{g}_0^{ss}}v,$$ where $v \in \mathbb{E}.$

{\it Proof.} See, e.g., Humphreys \cite{Humphreys}. $\Box$

Before we state the next lemma, let us do some comments on the
relationship between the Killing forms $(,)_{\mathfrak{g}_0^{ss}}$ and
$(,)_{\mathfrak{g}}$.
  It is well known that the restriction
of $(,)_{\mathfrak{g}}$ to $\mathfrak{g}_0^{ss}$ is a non degenerate
and obviously an invariant bilinear form, and therefore there is a
constant $\kappa \in \mathbb{C}^{\times},$ such that for $X,Y \in
\mathfrak{g}_0^{ss}$ we have $(X,Y)_{\mathfrak{g}_0^{ss}}=\kappa (X,
Y)_{\mathfrak{g}}$ - due to the uniqueness of invariant non-degenerate
forms up to a non zero complex multiple. The bases
$\{\eta^i\}_{i=1}^t$ and $\{\eta_i\}_{i=1}^t$ of
$\mathfrak{g}_0^{ss}$ are not dual with respect to the Killing form
$(,)_{\mathfrak{g}_0^{ss}}$ in general. For further purposes, we can consider these bases being also bases of the appropriate
complexified Lie algebras. According to the
relationship between  the Killing forms in question, we know that
$\{\eta^i\}_{i=1}^t$ and $\{\kappa^{-1}\eta_i\}_{i=1}^t$  are dual
with respect to $(,)_{\mathfrak{g}_0^{ss}}.$ We would like to
compute $(\sum_{i=1}^t \eta^i \eta_i).v.$ Due to Theorem 4, we can
write $(\sum_{i=1}^t\eta^i\kappa^{-1}\eta_i) .v = (\lambda, \lambda
+ 2 \delta)_{\mathfrak{g}_0^{ss}}v,$ if $v \in L(\lambda).$
Therefore $(\sum_{i=1}^t \eta^i
\eta_i).v=\kappa(\lambda,\lambda+2\delta)_{\mathfrak{g}_{0}^{ss}}v.$
Let us denote $(Gr,Gr)_{\mathfrak{g}}=:\rho^{-1},$ i.e.,
$\eta^{t+1}=Gr$ whereas $\eta_{t+1}= \rho Gr.$ Thus if  $Gr$ acts by
a complex number $c,$ we have that the action of $\eta^{t+1}\eta_{t+1}$
is by $\rho c^2.$
We will use these computations in the proof of the
following

{\bf Lemma 2:}
Suppose  $\mathbb{E}$ is an irreducible $\mathfrak{p}^{\mathbb{C}}$-module, the action of $(\mathfrak{g}_+)^{\mathbb{C}}$
being trivial and the highest weight of $\mathbb{E}$ over $(\mathfrak{g}_0)^{\mathbb{C}}$ is
$(\lambda,c) \in \mathfrak{h}^*\times \mathbb{C}.$
Let us further suppose that $\mathbb{E}\otimes (\mathfrak{g}_1)^{\mathbb{C}}$ decomposes into a finite direct sum
$\mathbb{E}\otimes \mathfrak{g}_1 = \bigoplus_{\mu}\mathbb{E}^{\mu}$ of irreducible $(\mathfrak{g}_0^{ss})^{\mathbb{C}}$-modules,  where $\mathbb{E}^{\mu}$ is an irreducible $(\mathfrak{g}_0^{ss})^{\mathbb{C}}$-module
with a highest weight $\mu.$ Let us fix a set of projections $\pi_{\mu}$ onto the irreducible summands in $\mathbb{E}\otimes (\mathfrak{g}_1)^{\mathbb{C}}.$
Assume further that $(\mathfrak{g}_1)^{\mathbb{C}}$ is an irreducible $(\mathfrak{g}_0^{ss})^{\mathbb{C}}$-module with a highest weight $\alpha.$
Then

\begin{equation}
\Psi  = \sum_{\mu} (\rho c-\kappa c_{\lambda \alpha}^{\mu})\pi_{\mu}
\end{equation}

{\it Proof.} Let us do the following computation with "Casimir"
operators \\ $\sum_{i=1}^{t+1} \eta^i\eta_i \in
\mathfrak{U}(\mathfrak{g}_0)$. For $X \in \mathfrak{g}_1$ and $v \in
\mathbb{E},$ we have:

\begin{eqnarray}
\nonumber\sum_{i=1}^{t+1}(\eta^i\eta_i).(X\otimes
v) &=&\sum_{i=1}^{t+1}(\eta^i\eta_i).X\otimes v+X\otimes
\sum_{i=1}^{t+1}(\eta^i\eta_i).v+\\
&&2\Psi(X\otimes v)\label{rov0},
\end{eqnarray}
where we have used  Lemma 1. Now, we would like to compute the first
two terms of the R.H.S. of the last written equation using the universal Casimir
element of $\mathfrak{g}_0^{ss}$, see Theorem  4.

\begin{eqnarray}
\sum_{i=1}^{t+1}(\eta^i\eta_i).X\otimes v
=\kappa(\alpha,\alpha+2\delta)_{\mathfrak{g}_0^{ss}}X\otimes v + \rho X\otimes v \label{rov1}
\end{eqnarray}

\begin{eqnarray}
X\otimes \sum_{i=1}^{t+1}(\eta^i\eta_i).v
=\kappa(\lambda,\lambda+2\delta)_{\mathfrak{g}_0^{ss}}X\otimes v+ \rho c^2X\otimes v \label{rov2}
\end{eqnarray}

Let us compute the L.H.S. of (\ref{rov0})

\begin{eqnarray}
\nonumber \sum_{i=1}^{t+1}(\eta^i\eta_i).(X\otimes
v)
&=&\sum_{\mu}\kappa(\mu, \mu+2\delta)_{\mathfrak{g}_0^{ss}}\pi_{\mu}(X\otimes
v)\\&&+\sum _{\mu}\pi_{\mu}[\rho X\otimes v+2\rho c X\otimes v+\rho c^2X\otimes
v]\label{rov3}
\end{eqnarray}

Substituting equations (\ref{rov1}), (\ref{rov2}) and (\ref{rov3})
into equation (\ref{rov0}) we obtain
\begin{eqnarray*}
\sum_{\mu}\kappa(\mu,\mu+2\delta)_{\mathfrak{g}_0^{ss}}\pi_{\mu}(X\otimes
v)+2\sum_{\mu}\rho c\pi_{\mu}(X\otimes v)+
\sum_{\mu}\rho c^2\pi_{\mu}(X\otimes v)+\rho X\otimes v =\\= 2\Psi(X\otimes
v)+\kappa(\alpha,\alpha+2\delta)_{\mathfrak{g}_0^{ss}}X\otimes v+\rho X\otimes
v+\kappa(\lambda,\lambda+2\delta)_{\mathfrak{g}_0^{ss}}X\otimes v+\rho c^2X\otimes v.&&
\end{eqnarray*}

As a result we obtain
$$\Psi(X\otimes v)=\sum_{\mu}(\rho c-\kappa c_{\lambda\alpha}^{\mu})\pi_{\mu}(X\otimes
v).$$

$\Box$

\subsection{Infinitesimal level classification}

Let $(\mathbb{V},\omega)$ be a real symplectic vector space of
dimension $2l+2,$ $l\geq 3$. In this subsection, we shall focus our
attention to the specific case of symplectic Lie algebra $
\mathfrak{sp}(\mathbb{V},\omega) \simeq
\mathfrak{sp}(2l+2,\mathbb{R})$    and its parabolic subalgebra
$\mathfrak{p}$ introduced in section 2.
 We shall be investigating the
vector space $\hbox{Hom}_{\mathfrak{p}}(J^1\mathbb{E},\mathbb{F})$ for suitable $\mathfrak{p}$-modules $\mathbb{E}, \mathbb{F},$
i.e., classify the  first order invariant differential operator at the infinitesimal level.
For a moment, we shall consider a complex setting.

The  complex
symplectic Lie algebra $\mathfrak{g}^{\mathbb{C}} =
\mathfrak{sp}(2l+2,\mathbb{C})$ possesses a $|2|$-grading,
$$\mathfrak{g}^{\mathbb{C}}  = \mathfrak{g}_{-2}^{\mathbb{C}} \oplus \mathfrak{g}_{-1}^{\mathbb{C}}\oplus
\mathfrak{g}_0^{\mathbb{C}}\oplus \mathfrak{g}_1^{\mathbb{C}} \oplus \mathfrak{g}_2^{\mathbb{C}},$$  such
that $\mathfrak{g}_{2}^{\mathbb{C}} \simeq \mathbb{C}$, $ \mathfrak{g}_1 ^{\mathbb{C}}\simeq
\mathbb{C}^{2l}$, $\mathfrak{g}_0^{\mathbb{C}} = (\mathfrak{g}_0^{ss})^{\mathbb{C}}\oplus
(\mathfrak{z}(\mathfrak{g}_0))^{\mathbb{C}} \simeq
\mathfrak{sp}(2l,\mathbb{C})\oplus \mathbb{C}.$ This splitting could
be displayed as follows.
Choose a basis $B$ of $\mathbb{V}$ such that $\omega,$ expressed in coordinates   with respect to $B,$
is given by $\omega((z^1,\ldots, z^{2l+2}),(w^1,\ldots, w^{2l+2}))=w^{1}z^{2l+2}+\ldots +w^{l+1}z^{l+2} - w^{l+2}z^{l+1}- \ldots - w^{2l+2}z^{1}.$
 For $A\in \mathfrak{sp}(2l+2,\mathbb{C})$
we have:

$$  A=\left(
\begin{tabular}{c|ccc|c}
$\mathfrak{g}_{0}$& &$\mathfrak{g}_{1}$&  & $\mathfrak{g}_2$ \\
\hline
&&&&\\
$\mathfrak{g}_{-1}$& &$\mathfrak{g}_0$& &$\mathfrak{g}_1$\\
&&&&\\
\hline $\mathfrak{g}_{-2}$& &$\mathfrak{g}_{-1}$& &$\mathfrak{g}_{0}$\\
\end{tabular}\right)$$ with respect to $B.$
%with respect to a basis $B$ with that property that the matrix
%of $\omega$ with respect to $B$ is an antidiagonal matrix having $1$'s in the upper-right side and $(-1)$'s at the lower-left part of the antidiagonal.
%has the form $\mathbb{J}=\left(
%\begin{tabular}{c|c}
%$0$& $1$ \\
%\hline
%$-1$ & $0$
% \\
%\end{tabular}\right)$
As one can easily compute, the parabolic subalgebra
$\mathfrak{p}^{\mathbb{C}}=(\mathfrak{g}_0)^{\mathbb{C}}\oplus
(\mathfrak{g}_1)^{\mathbb{C}} \oplus (\mathfrak{g}_2)^{\mathbb{C}}$
is a complexification of the Lie algebra of the group $P$ introduced
in section 2, where we have defined the metaplectic contact
projective geometry. Before we state the next theorem, we should
compute the coefficients $\rho$ and $\kappa$ for the   case
$\mathfrak{g}=\mathfrak{sp}(2l+2,\mathbb{C})$ considered with the
grading   given above. One can easily realize, that
$$Gr=\left(
\begin{tabular}{c|ccc|c}
1& & 0 &  & 0 \\
\hline
&&&&\\
0& &$0_{2l}$& &0\\
&&&&\\
\hline 0& & 0& &-1\\
\end{tabular}\right)$$
is the grading element, and that $(Gr,Gr)_{\mathfrak{g}}=8l+8.$
  Computing the square-norm of an element of $\mathfrak{g}_0^{ss}$ via $(,)_{\mathfrak{g}}$ and  $(,)_{\mathfrak{g}_0^{ss}},$ one obtains for the ratio $\kappa =\frac{l}{l+1}.$
    Further, let us introduce a bilinear form
     $<,>$ on $\mathfrak{h}^*,$ in which the orthogonal basis $\{\epsilon_i\}_{i=1}^l$ is orthonormal.
 The relation between the Killing form $(,)_{\mathfrak{g}_0^{ss}}$  and $<,>$ is given by
 $(X,Y)_{\mathfrak{g}_0^{ss}}=\frac{1}{8l+8}<X,Y>$ for   $X, Y \in \mathfrak{h}^*.$ For each $\lambda, \mu, \alpha \in \mathfrak{h}^*,$ let us define a complex number
  $$\tilde{c}_{\lambda \alpha}^{\mu}=\frac{1}{2}(<\mu,\mu+2\delta>-<\lambda,\lambda+2\delta>-<\alpha,\alpha + 2\delta>).$$
  Substituting the computed values  of $\rho$ and $\kappa$ and the relation between $(,)_{\mathfrak{g}_0^{ss}}$ and $<,>$ into formula (2), we obtain a prescription for mapping $\Psi$
  (in the metaplectic contact projective case)
  $$\Psi=\frac{1}{8l+8}\sum_{\mu}(c - \tilde{c}_{\lambda \alpha}^{\mu})\pi_{\mu}.$$

{\bf Theorem 5:} For $(\lambda, c), (\mu, d) \in A \times \mathbb{C},$
let $\mathbb{E}$ and $\mathbb{F}$ be two $\mathfrak{p}^{\mathbb{C}}$-modules such that $\mathbb{E}$ and $\mathbb{F}$
 are irreducible if considered as $(\mathfrak{g}_0)^{\mathbb{C}}$-modules with highest weight $(\lambda,c)$ and $(\mu,d),$
 respectively,
  and let $(\mathfrak{g}_+)^{\mathbb{C}}$ has a trivial action on each of these modules.
Further, suppose $\lambda \neq \mu.$
  Then
 $$\hbox{Hom}_{\mathfrak{p}^{\mathbb{C}}}(J^1\mathbb{E},\mathbb{F})  \simeq \left\{\begin{array}{l}
                                                       \mathbb{C} , \quad \hbox{if} \quad \mu \in
                                                       A_{\lambda} \quad \hbox{and} \quad d-1=c = \tilde{c}_{\lambda\varpi_1}^{\mu}
                                                                                                              \\
                                                       0 \quad {\hbox{in other cases.}}
                                                       \end{array}
                                                       \right.
                                                       $$

\noindent{\it Proof.}  Let us start with the second part of the statement, i.e.,
$\mu \notin A_{\lambda}$  or $c\neq \tilde{c}_{\lambda\varpi_1}^{\mu}$ or $d - 1\neq \tilde{c}_{\lambda\varpi_1}^{\mu},$ and consider
an element $T \in \hbox{Hom}_{\mathfrak{p}^{\mathbb{C}}}(J^1\mathbb{E},\mathbb{F}).$ Then
$T \in \hbox{Hom}_{(\mathfrak{g}^{ss}_0)^{\mathbb{C}}}(J^1\mathbb{E},\mathbb{F}).$
Because $T$ is a $\mathfrak{p}^{\mathbb{C}}$-homomorphism, we have that $T \in \hbox{Hom}_{(\mathfrak{g}_0^{ss} )^{\mathbb{C}}}(J^1_{R}\mathbb{E},\mathbb{F})$ due
to Theorem 3 (used in the complexified setting). We also know that
$$\hbox{Hom}_{(\mathfrak{g}_0^{ss} )^{\mathbb{C}}}(J^1_{R}\mathbb{E},\mathbb{F})= \hbox{Hom}_{(\mathfrak{g}_0^{ss} )^{\mathbb{C}}}(\mathbb{E}, \mathbb{F})\oplus\bigoplus_{\nu
\in A_{\lambda}}\hbox{Hom}_{(\mathfrak{g}_0^{ss}
)^{\mathbb{C}}}(L(\nu),L(\mu))$$ due to Theorem 2. If we suppose
$\mu \notin A_{\lambda}$ and $\lambda\neq \mu,$  then due to
Theorems 2.6.5, 2.6.6 in Dixmier \cite{Dixmier},  each member of the
direct sum is zero.  Now suppose that $\mu \in A_{\lambda}.$ Thus $c
\neq \tilde{c}_{\lambda \varpi_1}^{\mu}$ or $d-1 \neq
\tilde{c}_{\lambda \varpi_1}^{\mu}$. First suppose that $c \neq
\tilde{c}_{\lambda \varpi_1}^{\mu}.$ Using Theorem 2 and the cited
theorems of Dixmier, we see that $\hbox{Hom}_{(\mathfrak{g}_0^{ss}
)^{\mathbb{C}}}(J^1_{R}\mathbb{E},\mathbb{F}) \simeq
\hbox{Hom}_{(\mathfrak{g}_0^{ss} )^{\mathbb{C}}}(\mathbb{E},
\mathbb{F})\oplus\hbox{Hom}_{(\mathfrak{g}_0^{ss}
)^{\mathbb{C}}}(L(\mu),L(\mu))\simeq\hbox{Hom}_{(\mathfrak{g}_0^{ss}
)^{\mathbb{C}}}(L(\mu),L(\mu)),$ because the decomposition of
$(\mathfrak{g}_1)^{\mathbb{C}}\otimes \mathbb{E}$ is
multiplicity-free and $\lambda \neq \mu.$  Thus we can consider $T$
to be a $(\mathfrak{g}_0^{ss} )^{\mathbb{C}}$-intertwining operator
acting on the irreducible highest weight module $L(\mu)$. We have
two possibilities: $T: L(\mu) \to L(\mu)$ is either zero and we are
done, or $\hbox{Ker} \, T =\{0\}.$ We will suppose the latter
possibility.
 Take a nonzero
element $0\neq v \in L(\mu).$ Using the formula
 $\Psi = (8l+8)^{-1}\sum_{\nu}(c-\tilde{c}_{\lambda\varpi_1}^{\nu})\pi_{\nu},$ we obtain under the assumption
$c \neq \tilde{c}_{\lambda \varpi_1}^{\mu}$ that $\Psi(v)=(8l+8)^{-1}(c-\tilde{c}_{\lambda \alpha}^{\mu})v \neq 0.$  Because $\mbox{Ker} \,
T =\{0\},$ we have that $T\Psi(v) \neq 0$ and thus, according to
Theorem 3, $T$ it is not a $\mathfrak{p}^{\mathbb{C}}$-module homomorphism
because it does not vanish on the image of $\Psi.$ Secondly, consider
the case $d\neq \tilde{c}_{\lambda \alpha}^{\mu}+1.$ We can make the
following easy computation.  $d (S_1 \otimes v'')=Gr.(S_1\otimes
v'')=[Gr,S_1]\otimes v'' + S_1 \otimes Gr.v''=(1+c) S_1\otimes v''$
for $S_1 \in (\mathfrak{g}_1)^{\mathbb{C}}$ and $v'' \in \mathbb{E}.$ Thus $c=d-1$
and we are obtaining the case $c\neq \tilde{c}_{\lambda
\varpi_1}^{\mu},$ which was already handled.

Now, consider the case $\mu \in A_{\lambda},$ $c=\tilde{c}_{\lambda\alpha}^{\mu}$ and $d-1=\tilde{c}_{\lambda\alpha}^{\mu}$
and take a $T \in \hbox{Hom}_{\mathfrak{p}^{\mathbb{C}}}(J^1\mathbb{E},\mathbb{F}).$
As in the previous case, this implies $T \in \hbox{Hom}_{(\mathfrak{g}_0^{ss} )^{\mathbb{C}}}
(J^1_{R}\mathbb{E},\mathbb{F}).$
Decomposing  $J^1_{R}\mathbb{E}=L(\lambda) \oplus F(\varpi_1)\otimes L(\lambda)$ into irreducible modules  and
 substituting this decomposition into   $\hbox{Hom}_{(\mathfrak{g}_0^{ss} )^{\mathbb{C}}}(J^1_{R}\mathbb{E},\mathbb{F}),$
 we obtain a direct sum

$$\hbox{Hom}_{(\mathfrak{g}_0^{ss} )^{\mathbb{C}}}(\mathbb{E} ,\mathbb{F})\oplus\bigoplus_{\nu \in
A_{\lambda}}\hbox{Hom}_{(\mathfrak{g}_0^{ss} )^{\mathbb{C}}}(L(\nu),L(\mu)).$$

According to our assumptions $\mu \in A_{\lambda}$ and $\lambda \neq \mu,$ and due to the
structure of the set $A_{\lambda},$ we know that the direct sum
simplifies into a space isomorphic to $\mathbb{C}$ (using the above
cited theorem of Dixmier once more).
  Thus we know that
$\hbox{Hom}_{\mathfrak{p}^{\mathbb{C}}}(J^1\mathbb{E},\mathbb{F}) \subseteq
\hbox{Hom}_{(\mathfrak{g}_0^{ss} )^{\mathbb{C}}}(J^1_{R}\mathbb{E},\mathbb{F})\simeq
\mathbb{C}.$
 To obtain an equality in the previous inclusion, consider the one dimensional vector space of $(\mathfrak{g}^{ss}_0)^{\mathbb{C}}$-homomorphisms
 $\{w\widetilde{\pi_{\mu}}| w \in \mathbb{C}\},$
 where $\widetilde{\pi_{\mu}}$ is a trivial extension of the
 projection $(\mathfrak{g}_1)^{\mathbb{C}}\otimes \mathbb{E} \to L(\mu).$
The elements of   this vector space are clearly
$(\mathfrak{g}^{ss}_0)^{\mathbb{C}}$-homomorphisms,   which vanish
on the image of $\Psi,$  if $c=\tilde{c}_{\lambda\alpha}^{\mu},$ and
they factorize through the restricted jets. What remains is to show
that for each $w\in \mathbb{C},$ mappings $w\widetilde{\pi_{\mu}}$
are not only $(\mathfrak{g}_{0}^{ss})^{\mathbb{C}}$-homomorphisms,
but also $(\mathfrak{g}_0)^{\mathbb{C}}$-homomorphisms. Notice that
it is sufficient to test the condition only on
$(\mathfrak{g}_1)^{\mathbb{C}}\otimes \mathbb{E}$ because $Gr \in
(\mathfrak{g}_0)^{\mathbb{C}},$ and $\widetilde{\pi_{\mu}}$ is the
trivial extension, see formula (1). For $S_1\in
(\mathfrak{g}_1)^{\mathbb{C}}$ and $v'' \in \mathbb{E}$, we have
$Gr.\widetilde{\pi_{\mu}}(S_1\otimes v'')=d\widetilde{\pi_{\mu}}(S_1
\otimes v'')$ by definition.  Now, let us evaluate
$\widetilde{\pi_{\mu}}Gr.(S_1\otimes
v'')=\widetilde{\pi_{\mu}}([Gr,S_1]\otimes v'' + S_1 \otimes
Gr.v'')= \widetilde{\pi_{\mu}}(S_1 \otimes v'' + c S_1\otimes
v'')=(1+c)\widetilde{\pi_{\mu}}(S_1 \otimes v'')=d
\widetilde{\pi_{\mu}}(S_1\otimes
v'')=Gr.\widetilde{\pi_{\mu}}(S_1\otimes v''),$ thus
$\widetilde{\pi_{\mu}}$ commutes with the action of $Gr.$ Therefore
$\widetilde{\pi_{\mu}}$ is a
$(\mathfrak{g}_0)^{\mathbb{C}}$-homomorphism and   the statement
follows using Theorem 3.
 $\Box $

Let us remark, that  for $\lambda=\mu,$ the space of homomorphisms is also one dimensional. But this case leads to zeroth order operators, which   are
not interesting from the point of view of our classification.
Let us derive an easy corollary of the above theorem.

{\bf Corollary 1:}  The preceding theorem remains true for a real
form $\mathfrak{f}$ of $(\mathfrak{g}_0^{ss})^{\mathbb{C}},$ if one
considers complex representations and complex linear homomorphisms.
In particular, it remains true for the split real form
$\mathfrak{f}=\mathfrak{g}_0^{ss} \simeq \mathfrak{sp}(2l,
\mathbb{R}).$

{\it Proof.} First, observe that the decomposition of $F(\varpi_1)
\otimes L(\lambda)$ remains the same also over $\mathfrak{f}$. For it, let us take an irreducible summand $\mathbb{M}$ in the decomposition and suppose there is a proper nontrivial complex submodule
$\mathbb{M}'$ of   $\mathbb{M}.$  For $v \in
\mathbb{M}'$ and $X+iY\in \mathfrak{f}+i\mathfrak{f},$ we get that
$(X+iY).v=X.v+iY.v.$ Using the fact that $\mathbb{M}'$ is closed
under complex number multiplication and  $X.v, Y.v \in
\mathbb{M}',$   we would obtain that $\mathbb{M}'$ is
$(\mathfrak{g}_{0}^{ss})^{\mathbb{C}}$-invariant, which is a
contradiction.

Second, we would like to prove  that each
$\mathfrak{f}$-invariant complex linear endomorphisms
of an irreducible module, say $\mathbb{F},$ is a scalar. It is easy
to observe, that such an endomorphism is actually
$(\mathfrak{g}_0^{ss})^{\mathbb{C}}$-endomorphism, i.e., the theorem of Dixmier
used in the proof of the previous theorem, could be applied  and the
corollary follows. $\Box$

\subsection{Globalized level classification}

In this subsection, we shall extend the results obtained in the previous one to the group level.
We will do it using some basic facts on globalization techniques.

Let $(\mathbb{V},\omega)$ be a real symplectic vector space of
real dimension $2l+2,$ $l\geq 3,$ $G=Sp(\mathbb{V},\omega)$ and $P$ as
described in section 2. First, we   introduce the groups, we
shall be considering. Let $G_+, G_0, G_0^{ss},K$ be the unipotent
part, the Levi factor, the semisimple part of $P$ and the maximal
compact subgroup of $G,$ respectively. Recall that we have fixed a
non-trivial 2-fold covering $q: \tilde{G} \to G$ of the symplectic
group $G$ by the metaplectic group
$\tilde{G}=Mp(\mathbb{V},\omega).$ Let us denote the respective
$q$-preimages by $\tilde{G}_+, \tilde{G}_0, \tilde{G}_0^{ss},
\tilde{K}.$ Further, let us denote the maximal compact subgroup of
the semisimple part $G_0^{ss}$ of the Levi factor by $K_0^{ss}$ and
its $q$-preimage by $\tilde{K}_0^{ss}.$ We have $$\tilde{K}_0^{ss}
\simeq \widetilde{U(l)}=\{(u,z) \in U(l) \times
\mathbb{C}^{\times}| \hbox{det} \, u=z^2\},$$ which is obviously
connected, see Tirao, Vogan and Wolf \cite{TVW}.

Second, let us introduce a class of $\tilde{P}$-modules we shall be dealing with.
In Kashiwara, Vergne \cite{KV}, the so called metaplectic (or Segal-Shale-Weil or oscillator) representation over $\tilde{G}_0^{ss}$ is introduced.
Let ${\bf S_+}$ be the irreducible submodule of the Segal-Shale-Weil representation consisting of even functions.
 Let us take the underlying $(\mathfrak{g}_{0}^{ss},\tilde{K}_0^{ss})$-module and denote it by $S_+.$ The $\mathfrak{g}_0^{ss}$-module structure of this representation coincides with the irreducible highest weight module structure of $\mathbb{S}_+,$ which   was introduced in section 3.
 For a choice of a weight $\lambda \in A,$ we know that there exists
a dominant integral weight $\nu$ (with respect to choices made in
section 3), such that $\mathbb{L}:=L(\lambda) \subseteq
\mathbb{S}_+\otimes F(\nu).$ Because $\mathbb{S}_+\otimes F(\nu)$
decomposes without multiplicities, we have an identification of
$L(\lambda)$ with its isomorphic module in $\mathbb{S}_+\otimes
F(\nu).$  Now we would like to make $\mathbb{L}$ a
$(\mathfrak{g}_0^{ss},\tilde{K}_0^{ss})$-module. Using a result of Baldoni \cite{WB}, this could be
done as follows.  Because $S_+$ and $F(\nu)$ are $(\mathfrak{g}_0^{ss}, \tilde{K}_0^{ss})$-modules, their tensor product
is a $(\mathfrak{g}_0^{ss}, \tilde{K}_0^{ss})$-module as well.
Using the fact that $\tilde{K}_0^{ss}=\widetilde{U(l)}$ is connected, we are obtaining a $(\mathfrak{g}_0^{ss}, \tilde{K}_0^{ss})$-module  structure
on each irreducible summand in $\mathbb{S}_+\otimes F(\nu),$ in particular on $\mathbb{L}.$
 Denote the resulting $(\mathfrak{g}_0^{ss},\tilde{K}_0^{ss})$-module by $L.$
 Using   globalization results of Kashiwara and Schmid in \cite{KS},
 there exists a  minimal globalization for this $(\mathfrak{g}_0^{ss},\tilde{K}_0^{ss})$-module,
 which will be denoted by ${\bf L}=:{\bf L}(\lambda).$ (For this topic, see also   Vogan \cite{Vogan}
and Schmid \cite{Schmied}.)
 Thus ${\bf L}(\lambda)$ is a complex $\tilde{G}_0^{ss}$-module. Further, we need to specify the action of the center of $\tilde{G}_0$ and that one of the unipotent part $\tilde{G}_+.$
For each $(\lambda, c) \in A \times \mathbb{C}$ we suppose, that the
unipotent $\tilde{G}_+$ acts trivially on ${\bf L}(\lambda)$ and the
grading element $Gr$ in  the Lie algebra of the center of the Levi
factor $\tilde{G}_0$ acts by multiplication by a complex number $c
\in \mathbb{C}.$  Since the center is isomorphic to $\mathbb{R}^{\times}$ we need to specify
the action of, e.g., $-1 \in \mathbb{R}^{\times}.$ This action should be any  $\gamma \in \mathbb{R}$ satisfying $\gamma^2=1.$
So we have obtained a $\tilde{P}$-module structure
on ${\bf L}(\lambda)$ which we will refer to as ${\bf
L}(\lambda,c)_{\gamma}.$
% This module is an admissible  module over the identity component of
%$\tilde{P}$. Sometimes it is more convenient to introduce projectivizations of the groups of questions %instead of speaking of
%"modules over the identity component", because the resulting projectivized parabolic group is %connected. Moreover, taking the projectivizations makes
%the action of $G$ on $\mathbb{PV}$ effective.
Let
us remark, that defining the action of $\tilde{G}_+$ to be  trivial,
is actually no restriction, when one considers only irreducible
admissible $\tilde{P}$-modules. We shall call the corresponding
associated bundles  higher symplectic bundles and the corresponding
$1^{st}$ order invariant differential operators  symplectic spinor
operators, stressing the fact that the representations of $\tilde{P}$
we are considering are coming from   higher symplectic spinor
modules.

{\bf Theorem 6.} Let $(\lambda,c,\gamma), (\mu,d, \gamma') \in A\times \mathbb{C} \times \mathbb{Z}_2$
\footnote{The group $\mathbb{Z}_2$ is considered as multiplicative, i.e., $\mathbb{Z}_2=\{-1,1\}.$},
$\lambda \neq \mu$ and $(p:\tilde{\mathcal{G}} \to M^{2l+1}, \omega)$ be a metaplectic contact projective geometry of rank $l$.
Consider the $\tilde{P}$-modules ${\bf E}:={\bf L}(\lambda, c)_{\gamma}$ and
${\bf F}:={\bf L}(\mu,d)_{\gamma'}.$
  Then for the vector space of invariant differential operators up to a zeroth order we have

$$\hbox{Diff}({\bf E}M,{\bf F}M)_{(p:\tilde{\mathcal{G}}\to M^{2l+1},  \omega)}^1  \simeq \left\{\begin{array}{l}
                                                       \mathbb{C} \quad  \hbox{if} \  \mu \in
                                                       A_{\lambda} \hbox{,} \
                                                       d-1=c = \tilde{c}_{\lambda\varpi_1}^{\mu} \ \hbox{and} \ \gamma=\gamma'\\
                                                       0 \quad {\hbox{in other cases.}}
                                                       \end{array}
                                                       \right.
                                                       $$

{\it Proof.} According to the definition of first order invariant differential operators between sections of associated vector bundles over Cartan geometries,
the vector space $\hbox{Diff}({\bf E}M,{\bf F}M)_{(p:\tilde{\mathcal{G}}\to M,\omega)}^1$ is isomorphic to the space $\hbox{Hom}_{\tilde{P}}(J^1{\bf E},{\bf F}).$
From the definition of the minimal globalization, it follows that it gives a natural bijection between Hom's of respective categories: the Harish-Chandra category of
$(\mathfrak{p},\widetilde{K \cap G_0})$-modules and the category of admissible $\tilde{P}$-modules, see
Kashiwara, Schmid \cite{KS}. Thus we have
$\hbox{Hom}_{\widetilde{P}}(J^1{\bf E},{\bf F}) \simeq \hbox{Hom}_{(\mathfrak{p},\widetilde{K\cap G_0})}(J^1E,F).$
Because the identity component $(\widetilde{K\cap G_0})_1$ is connected by definition, we can write
$\hbox{Hom}_{(\mathfrak{p},(\widetilde{K\cap G_0})_1)}(J^1 E, F)\simeq \hbox{Hom}_{\mathfrak{p}}(J^1{\mathbb E},\mathbb{F}),$ see W. Baldoni \cite{WB}.
 It remains to show that each $\mathfrak{p}$-module homomorphism is actually
a $(\mathfrak{p},(\widetilde{K\cap G_0})_{-1})$-module homomorphism,
where $(\widetilde{K\cap G_0})_{-1}$ denotes the component of the
group $\widetilde{K\cap G_0}$ to which $-1$ belongs. Let us
parameterize the elements of the $(-1)$-component of
$(\widetilde{K\cap G_0}) \simeq \widetilde{U(l)}\times \mathbb{Z}_2$
by  pairs $(k,-1),$ $k\in \widetilde{U(l)},$ and denote the
appropriate $\tilde{P}$-representation on $E$ by $\rho.$ We can
easily check that for $(v',S \otimes v'')\in J^1{\bf E},$ we have
$(k,-1).(v',S\otimes v'')=(\rho(k,-1)v',Ad(k,-1)S\otimes \rho(k,-1)
v'')=(\gamma\rho(k,1)v',Ad(k,1)S\otimes \gamma\rho(k,1)v'')=\gamma
(k,1).(v',S\otimes v'').$ Further for a $\mathfrak{p}$-homomorphism
$T \in \hbox{Hom}_{\mathfrak{p}}(J^1{\mathbb E}, {\mathbb F}),$ we
can write $T(k,-1).(v',S\otimes v'')=\gamma T(k,1).(v',S \otimes
v'')=\gamma (k,1).T(v',S\otimes v'')=\gamma
\gamma'(k,-1).T(v',S\otimes v'').$ Thus we have also
$\hbox{Hom}_{(\mathfrak{p},\widetilde{K\cap G_0})}(J^1 E, F)\simeq
\hbox{Hom}_{\mathfrak{p}}(J^1{\mathbb E},\mathbb{F})$ if
$\gamma=\gamma'.$ The Hom at the right hand side was determined in
corollary 1. In the case $\gamma \neq \gamma',$ we have that $T=0$
and the proof is finished. $\Box$

\section{ Examples: contact projective Dirac, twistor and Rarita-Schwinger operators}

 In this section, we shall introduce three main examples of contact projective analogues of   Dirac, twistor and Rarita-Schwinger operators known from Riemannian  and partly from symplectic geometry. In each of the next paragraphs, we suppose that a metaplectic contact projective geometry
  $(p:\tilde{\mathcal{G}} \to M^{2l+1}, \omega)$ of rank $l$ is fixed.

\paragraph{Contact projective Dirac operator.}
For $\lambda = -\frac{1}{2}\varpi_l,$
 we have $A_{\lambda}=\{\varpi_1-\frac{1}{2}\varpi_l, \varpi_{l-1} - \frac{3}{2}\varpi_l\}$ according to Theorem 2.
Take  $\mu=\varpi_{l-1}-\frac{3}{2}\varpi_l  \in A_{\lambda}.$ Using
$\delta= l\epsilon_1 + (l-1)\epsilon_2 + \ldots + \epsilon_l,$  we
obtain that $\tilde{c}_{\lambda \alpha}^{\mu} =\frac{1+2l}{2}.$ Thus
for conformal weight $c=\frac{1+2l}{2}$ and $\gamma \in \mathbb{Z}_2$ there is an invariant
differential operator $\mathfrak{D}^{\frac{1}{2}}:\Gamma(M,{\bf
L}(\lambda,\frac{1+2l}{2})_{\gamma}M) \to \Gamma(M^{2l+1},{\bf
L}(\mu,\frac{3+2l}{2})_{\gamma}M).$ This operator could be called  contact
projective   Dirac operator  because of the analogy with the
orthogonal case.

%Vysledky:
%$j: \mathfrak{g}_0^{ss} \to \mathfrak{g}_0$
%$\rho: \mathfrak{p} \to \hbox{End}(L(\lambda))$
%$z^i$ $z_i$ baze $\mathfrak{g}_0^{ss}$ dualni vuci $(,)_m$
%$j(z_i)$ $j(z^i)$ vektory v $\mathfrak{g}$ ale nejsou dualni vuci $(,)_{v=velka}$
%$j(z^i)j(z_i)v=\frac{l}{l+1}(\lambda,\lambda + 2 \delta)_{Kill, m=mala}$
%VIM:$(8l+8)(,)_m=8l(,)_v$
%$(j(z^i),j(z_i)\frac{l+1}{l})_m=\frac{l+1}{l}(j(z^i),j(z_i))_m=(j(z^i),j(z_i))_v=\delta_{ij}.$ Takto dostavam bazi dualni vuci
%$(,)_m$ a pro ni vim, ze
%$j(z^i) j(z_i) \frac{l+1}{l} = (\lambda,\lambda+2\delta)$
%Ale ja chci $j(z^i) j(z_i)=C.$ Odtud
%$C=\frac{l}{l+1}(\lambda,\lambda+2\delta).$

\paragraph{Contact projective twistor operator.}
Taking the same $\lambda=-\frac{1}{2}\varpi_l$ as in the previous
example and $\mu = \varpi_1-\frac{1}{2}\varpi_l,$ we obtain
$c=\frac{1}{2}$ and the corresponding operator
$\mathfrak{T}:\Gamma(M,{\bf L}(\lambda,\frac{1}{2})_{\gamma}M) \to \Gamma(M,
{\bf L}(\mu,\frac{3}{2})_{\gamma}M)$ ($\gamma \in \mathbb{Z}_2$) is called contact projective  twistor
operator also due to the analogy with the orthogonal case.

\paragraph{Contact projective Rarita-Schwinger operator.}
Here, take $\lambda = \varpi_1 - \frac{1}{2}\varpi_l.$ $A_{\lambda}=\{\varpi_2-\frac{1}{2}\varpi_l, 2\varpi_1-\frac{1}{2}\varpi_l, -\frac{1}{2}\varpi_l, \varpi_1+\varpi_{l-1}-\frac{3}{2}\varpi_{l}\}.$ For $\mu = \varpi_1+\varpi_{l-1}-\frac{3}{2}\varpi_l,$ we obtain
$c =\frac{1+2l}{2},$ and we shall call this operator contact projective Rarita-Schwinger operator,\\
$\mathfrak{D}^{\frac{3}{2}}:\Gamma(M,{\bf L}(\lambda,\frac{1+2l}{2})_{\gamma}M) \to \Gamma(M,{\bf L}(\mu,\frac{3+2l}{2})_{\gamma}M),$ where again $\gamma \in \mathbb{Z}_2.$
\\

{\bf Remark:} It   may be interesting to mention, that computing formally the conformal weights using a Lepowsky generalization of a result of Bernstein-Gelfand-Gelfand on homomorphism of non-true Verma-modules, one gets exactly the same weights, although
  Lepowsky is considering only Verma modules induced by  finite dimensional representations.

\end{document}